\documentclass[11pt]{article}
\usepackage[top=1in,bottom=1in,left=1 in,right=.8in,footskip=0.4in]{geometry}
\usepackage{multicol}

\usepackage{amsmath, amsthm, amssymb,mathtools}

\usepackage{authblk}
\usepackage[mathscr]{eucal}

\usepackage{xfrac}

\usepackage{color,soul,enumerate}
\usepackage{hyperref}

\usepackage{orcidlink}

\usepackage{tikz}
\usepackage{tikz-cd}
\usetikzlibrary{graphs}
\usetikzlibrary{graphs.standard}
\usetikzlibrary{arrows.meta}
\usepackage{calc}
\usepackage{float}
\usetikzlibrary{calc}
\tikzstyle{vertex}=[circle, draw, fill=black, inner sep=0pt, minimum size=3pt]
\tikzstyle{grvert}=[circle, draw, gray, fill=gray,  inner sep=0pt, minimum size=6pt]
\newcommand{\vertex}{\node[vertex]}

\newtheorem{theorem}{Theorem}

\newtheorem{obs}[theorem]{Observation}

\theoremstyle{definition}
\newtheorem{definition}[theorem]{Definition}

\theoremstyle{remark}

\newtheorem{counterx}[theorem]{Counterexample}
\newtheorem{claim}[theorem]{Claim}

\def\N{\mathcal{N}}
\def\G{\mathcal{G}}
\def\K{\mathcal{K}}
\def\V{\mathscr{V}}

\numberwithin{theorem}{section}
\numberwithin{equation}{section}
\theoremstyle{plain}

\begin{document}

	\title{A note on an application of discrete Morse theoretic techniques on the complex of disconnected graphs}
	\author[a]{Anupam Mondal~\orcidlink{0000-0002-6547-4835}\thanks{\texttt{anupam.mondal@tcgcrest.org}}}
	
	\author[b]{Pritam Chandra Pramanik\thanks{\texttt{pritam.pramanik.80@tcgcrest.org}}}
	
	\affil[a,b]{\small Institute for Advancing Intelligence (IAI), TCG CREST, Kolkata -- 700091, West Bengal, India}
	\affil[a]{\small Academy of Scientific \& Innovative Research (AcSIR), Ghaziabad -- 201002, Uttar Pradesh, India}
	\affil[b]{\small Department of Mathematics, National Institute of Technology (NIT) Durgapur, Durgapur -- 713209, West Bengal, India}
	\date{}
	\maketitle
	
	\begin{abstract}	  
	Robin Forman's highly influential 2002 paper \textit{A User's Guide to Discrete Morse Theory} presents an overview of the subject in a very readable manner. As a proof of concept, the author determines the topology (homotopy type) of the abstract simplicial complex of disconnected graphs of order $n$ (which was previously done by Victor Vassiliev using classical topological methods) using discrete Morse theoretic techniques, which are purely combinatorial in nature. The techniques involve the construction (and verification) of a discrete gradient vector field on the complex. However, the verification part relies on a claim that doesn't seem to hold. In this note, we provide a couple of counterexamples against this specific claim. We also provide an alternative proof of the bigger claim that the constructed discrete vector field is indeed a gradient vector field. Our proof technique relies on a key observation which is not specific to the problem at hand, and thus is applicable while verifying a constructed discrete vector field is a gradient one in general.
	
	\end{abstract}
		
	\section{Introduction}
	Topological properties of the abstract simplicial complex of disconnected graphs of order $n$, denoted by $\N_n$, became a focus of intense interest due to their relation to knot theory. As explained in~\cite{vassiliev}, the complex $\N_n$~\footnote{rather, to be more precise, the quotient complex $\sfrac{\G_n}{\N_n}$ (called the complex of connected graphs), where $\G_n$ is the complex of all graphs of order $n$} is an essential part of Vassiliev's spectral sequence calculation of the homology of the space of singular knots. The following theorem determining the homotopy type of $\N_n$ is also due to Vassiliev, who proved it using classical topological techniques.
	\begin{theorem}\cite{forman02,vassiliev}\label{thm}
		The complex $\N_n$ of disconnected graphs of order $n$ is homotopy equivalent to the wedge of $(n - 1)!$ spheres of dimension $(n - 3)$.
	\end{theorem}
	
	Robin Forman developed discrete Morse theory as a combinatorial analogue of (smooth) Morse theory in~\cite{forman}, which turned out to be immensely useful over the years. Given a finite (abstract) simplicial complex (or a finite regular CW complex), discrete Morse theory helps us capture the topology of the complex through an ``efficient" CW decomposition (i.e., one with many fewer cells than in the original decomposition) for the complex. A later (expository) article by Forman titled \textit{A User's Guide to Discrete Morse Theory}~\cite{forman02} presents an overview of the subject in a comprehensible and engaging manner. In the preface of the textbook \emph{Discrete Morse Theory}~\cite{scoville} (Student Mathematical Library, AMS), the author, Nicholas Scoville refers to this article as ``(still) the gold standard in the field." As a proof of concept, in Section~5 of the User's Guide$\ldots$, the author proves Theorem~\ref{thm} using only discrete Morse theoretic techniques, which are purely combinatorial in nature. The techniques involve the construction of a suitable (discrete) gradient vector field $\V$ on $\N_n$. However, a claim (see Claim~\ref{claim-forman}) is made while justifying that the constructed discrete vector field $\V$ is indeed a gradient vector field on $\N_n$, which doesn't seem to hold. In Section~\ref{main} of this note, we provide a couple of counterexamples that show the invalidity of the claim in two different ways.
	
	Here we remark that Dmitry Kozlov provides a brief justification of $\V$ being a gradient vector field in the book \emph{Organized Collapse: An Introduction to Discrete Morse Theory} (Graduate Studies in Mathematics, AMS)~\cite[Chapter~11]{kozlov}. However, we also provide an alternative justification in Section~\ref{alt-proof}. Our proof technique relies on an observation (Observation~\ref{obs}) mentioned in Section~\ref{concl} which is not specific to the problem at hand. Thus this technique may be adapted while verifying a constructed discrete vector field is a gradient vector field in various other setups in general. For example, it is used in past joint work~\cite{cmms,mms} by the first author of the present article.
	
	\section{Preliminaries}
	In this note, by a \emph{graph} (\emph{of order $n$}), we always mean a simple, finite, undirected, labelled graph on the vertex set $\{1,\ldots,n\}$. The edge set of a graph $G$ is a collection of $2$-subsets of its vertex set, and we denote it by $E(G)$. Since, we consider labelled graphs, two graphs $G_1$ and $G_2$ of the same order are the same if and only if $E(G_1)=E(G_2)$. In other words, we consider a pair of graphs to be different if they (i.e., their edge sets) are different up to labelling, even if they are the same \emph{up to a graph isomorphism}.
	We denote the edge consisting of the pair of vertices $i$ and $j$, i.e., the edge \emph{joining} $i$ and $j$, by $(i,j)$ (or, $(j,i)$). If $(i,j)$ is an edge in $G$, then by $G-(i,j)$, we mean the graph (of the same order) with the edge set $E(G)\setminus \{(i,j)\}$. On the other hand, if $(i,j) \notin E(G)$, then $G+(i,j)$ denotes the graph (of the same order) with the edge set $E(G)\cup \{(i,j)\}$.
	
	An (\emph{abstract simplicial}) \emph{complex} is a (finite, nonempty) collection, say $\K$, of finite sets with the property that every subset (including the empty set) of a set in the collection is also in the collection $\K$. A member of the collection $\K$ is called a \emph{simplex}  of $\K$. The \emph{dimension} of a simplex $\sigma$, denoted by $\dim(\sigma)$, is the number $|\sigma|-1$. If $\dim(\sigma)=d$, then $\sigma$ is a \emph{$d$-dimensional simplex} (or simply, a \emph{$d$-simplex}). We denote a $d$-simplex $\sigma$ by $\sigma^{(d)}$ whenever necessary.
	
	We note that the collection of the edge sets of all disconnected (labelled) graphs on the vertex set $\{1,\ldots,n\}$ is an abstract simplicial complex, and we denote it by $\N_n$. Abusing the notation, we consider the disconnected graphs themselves (instead of their edge sets) as the simplices of $\N_n$.
	
	\subsection{The basics of discrete Morse theory}
	First, we introduce the notion of a \emph{discrete vector field} and a (discrete) \emph{gradient vector field} on a complex as described in~\cite{forman,forman02}. 
	\begin{definition}[Discrete vector field]
		A discrete vector field $\V$ on a complex $\K$ is a collection of ordered pairs of simplices of the form $(\alpha,\beta)$ such that 
		\begin{enumerate}[(i)]
			\item $\alpha \subsetneq \beta$,
			\item $\dim(\beta)=\dim(\alpha)+1$,
			\item each simplex of $\K$ is in \emph{at most} one pair of $\V$.
		\end{enumerate}
	\end{definition}
	If the simplex $\alpha^{(p)}$ is paired off with the simplex $\beta^{(p+1)}$ in $\V$ (i.e., $(\alpha, \beta) \in \V$), then sometimes we pictorially represent it as $\alpha \rightarrowtail \beta$, and we say $\alpha$ (or $\beta$) is the \emph{tail} (or \emph{head}) \emph{of an arrow} in $\V$.
	
	Given a discrete vector field $\V$ on a complex $\K$, a $\V$-path is a sequence of simplices 
	\[\alpha_0^{(d)}, \beta_0^{(d+1)}, \alpha_1^{(d)}, \beta_1^{(d+1)}, \ldots, \alpha_k^{(d)}, \beta_k^{(d+1)}, \alpha_{k+1}^{(d)}\]
	such that for each $i \in \{0,\ldots,k\}$, the pair $(\alpha_i,\beta_i) \in \V$ and $\beta_i \supsetneq \alpha_{i+1} \ne \alpha_i$. Such a path is a \emph{non-trivial closed path} if $k \ge 0$ and $\alpha_{k+1} = \alpha_0$.
	\begin{definition}[Gradient vector field]
		A gradient vector field on a complex $\K$ is a discrete vector field $\V$ on $\K$ which does not admit non-trivial closed $\V$-paths.
	\end{definition}
		
		\begin{definition}[Critical simplex]
			Let $\V$ be a gradient vector field on a complex $\K$. A nonempty simplex $\alpha$ is a \emph{critical simplex} (with respect to $\V$) if one of the following holds:
	\begin{enumerate}[(i)]
		\item $\alpha$ does not appear in any pair of $\V$, or
		\item $\alpha$ is a $0$-simplex and $(\emptyset,\alpha) \in \V$.
	\end{enumerate}
		\end{definition}
	We note that gradient vector fields can also be realized as, and often referred to in literature as \emph{acyclic matchings}.

	We recall that a CW complex is a topological space built recursively by gluing cells (which are homeomorphic copies of balls) of increasing dimension. The fundamental theorem of discrete Morse theory~\cite{forman,forman02} states that if $\K$ is a complex and $\V$ is a gradient vector field on $\K$, then (the geometric realization of) $\K$ is homotopy equivalent to a CW complex with exactly one cell of dimension $p$ for each critical simplex (with respect to $\V$) of dimension $p$. The following is a useful corollary of the fundamental theorem.
	\begin{theorem}\cite{forman02}\label{wedge}
		If $\K$ is a complex and $\V$ is a gradient vector field on $\K$ such that the only critical simplices are one $0$-simplex and $k$ simplices of dimension $d$, then $\K$ is homotopy equivalent to the wedge of $k$ spheres of dimension $d$.
	\end{theorem}

	\section{Construction of a discrete vector field $\V$ on $\N_n$}
	For the sake of completeness, we describe the construction of a discrete vector field $\V$ on $\N_n$, as done in~\cite[Section~5]{forman02}, in brief here. Construction of $\V$ is done in stages, by first considering the edge $(1,2)$ and then the vertices $3, 4,  \ldots$ in the increasing order. At the first stage we have the following discrete vector field:
	\[\V_{12}=\{(G,G+(1,2)): G \in \N_n, (1,2) \notin E(G), G+(1,2) \in \N_n\}.\]
	We note that any $G \in \N_n$, which is unpaired in $\V_{12}$, has two (connected) components with the vertices $1$ and $2$ belonging to different components, and vice versa.
	
	Next, we consider a graph $G$ that is unpaired in $\V_{12}$ such that the vertices $1$ (or $2$) and $3$ are in the same component of $G$, and $G$ doesn't contain the edge $(1,3)$ (or $(2,3)$). We extend $\V_{12}$ by pairing off all such $G$ with $G + (1,3)$ (or $G + (2,3)$). Let $\V_3$ denote the resulting discrete vector field. We note that any $G \in \N_n$, which is unpaired in $\V_3$, contains the edge $(1,3)$ (or $(2,3)$), and have the property that the graph $G-(1,3)$ (or $G-(2,3)$) has three components, one containing the vertex $1$, one containing the vertex $2$, and one containing the vertex $3$ (and vice versa).
	
	Next, we consider a graph $G$ that is unpaired in $\V_3$, and pair it with $G + (1,4)$, or $G + (2, 4)$, or $G + (3, 4)$ if possible (we note that \emph{at most} one of these graphs is unpaired in $\V_3$). Let $\V_4$ denote the resulting extended discrete vector field. We continue in this fashion, considering the vertices $5$, $6,\ldots,n$ in the increasing order. Let $\V_i$ denote the discrete vector field constructed after the consideration of vertex $i$, and finally $\V \coloneqq \V_n$.
	
	\subsection{Proof of Theorem~\ref{thm} assuming $\V$ to be a gradient vector field}
	As demonstrated in~\cite[Section 5]{forman02}, with the assumption that $\V$ is a gradient vector field on $\N_n$, we may determine the homotopy type of $\N_n$ and settle Theorem~\ref{thm}. First, we may verify that the only unpaired graphs in $\V$ are the forests with two components (i.e., a disjoint union of two trees), with the vertices $1$ and $2$ belonging to different components, and both component trees have the property that the vertex labels are increasing along every ray starting from the vertex $1$ or the vertex $2$. There are exactly $(n - 1)!$ such graphs, and they each have ($n - 2$) edges. Thus, there are $(n-1)!$ critical $(n-3)$-simplices, only one critical $0$-simplex (viz., the graph with the single edge $(1,2)$), and no other critical simplices. Consequently, Theorem~\ref{thm} follows from Theorem~\ref{wedge}.
	
	\section{Counterexamples against a pivotal claim used while proving $\V$ is a gradient vector field}\label{main}
	While justifying that $\V$ is a gradient vector field, Forman first verifies that $\V_{12}$ is a gradient vector field. This follows from the observation that if $\alpha_0, \beta_0, \alpha_1$ is a $\V_{12}$-path, then $\alpha_1$ is the head of an arrow in $\V_{12}$, and thus the $\V_{12}$-path cannot be continued further. Next, an analogous argument for $\V$ is provided, which entails the following.
	\begin{claim}[{\cite[pp.~20--21]{forman02}}]\label{claim-forman}
		Let $\gamma = \alpha_0, \beta_0, \alpha_1$ denote a $\V$-path. In particular, $\alpha_0$ and $\beta_0$ must be paired in $\V$. The reader can check that if $\alpha_0$ and $\beta_0$ are first paired in $\V_i$, $i \ge 3$, then either $\alpha_1$ is the head of an arrow in $\V_i$, in which case the $\V$-path cannot be continued, or $\alpha_1$ is paired in $\V_{i-1}$. It follows by induction that there can be no closed $\V$-paths.
	\end{claim}
	However, this claim doesn't seem to hold. We provide a couple of counterexamples below.
	\begin{counterx}
		It may very well happen that $\alpha_0$ and $\beta_0$ are first paired in $\V_i$ for some $i \ge 3$, but $\alpha_1$, which is not the head of an arrow, gets first paired in $\V_j$ with $j>i$. A specific example of this happening for $n=5$ is as below in Figure~\ref{counter}.
\begin{figure}[!ht]
	\centering
	\begin{tikzpicture}
		\begin{scope}
			\draw[gray!50](0,0.7) circle (0.2);
			\draw [rounded corners,gray!50] (0,0.3)--(-0.7,-0.7)--(0,-1.7)--(0.7,-0.7)--cycle;
			
			\node at (0,-2) {$\alpha_0$};
			
			\vertex (v1) at (0,0.7) {};\node at (0.3,0.7) {$1$};
			\vertex (v2) at (0,0) {};\node at (0.3,0) {$2$};
			\vertex (v3) at (-0.5,-0.7) {};\node at (-0.8,-0.7) {$3$};
			\vertex (v4) at (0.5,-0.7) {};\node at (0.8,-0.7) {$4$};
			\vertex (v5) at (0,-1.4) {};\node at (0.3,-1.4) {$5$};
			
			\path
			(v2) edge (v4)
			(v4) edge (v5)
			(v5) edge (v3)
			;
		\end{scope}
	\draw [to reversed-to](1,-0.3) -- (2,-0.3);
	
	\begin{scope}[xshift=3cm]
		\draw[gray!50](0,0.7) circle (0.2);
		\draw [rounded corners,gray!50] (0,0.3)--(-0.7,-0.7)--(0,-1.7)--(0.7,-0.7)--cycle;
		
		\node at (0,-2) {$\beta_0$};
		
		\vertex (v1) at (0,0.7) {};\node at (0.3,0.7) {$1$};
		\vertex (v2) at (0,0) {};\node at (0.3,0) {$2$};
		\vertex (v3) at (-0.5,-0.7) {};\node at (-0.8,-0.7) {$3$};
		\vertex (v4) at (0.5,-0.7) {};\node at (0.8,-0.7) {$4$};
		\vertex (v5) at (0,-1.4) {};\node at (0.3,-1.4) {$5$};
		
		\path
		(v2) edge (v4)
		(v4) edge (v5)
		(v5) edge (v3)
		(v3) edge (v2)
		;
	\end{scope}

	\begin{scope}[xshift=6cm]
		\draw[gray!50](0,0.7) circle (0.2);
		\draw[gray!50](0,0) circle (0.2);
		\draw [rounded corners,gray!50] (-0.8,-0.55)--(0.8,-0.55)--(0,-1.7)--cycle;
		
		\node at (0,-2) {$\alpha_1$};
		
		\vertex (v1) at (0,0.7) {};\node at (0.3,0.7) {$1$};
		\vertex (v2) at (0,0) {};\node at (0.3,0) {$2$};
		\vertex (v3) at (-0.5,-0.7) {};\node at (-0.8,-0.7) {$3$};
		\vertex (v4) at (0.5,-0.7) {};\node at (0.8,-0.7) {$4$};
		\vertex (v5) at (0,-1.4) {};\node at (0.3,-1.4) {$5$};
		
		\path
		(v4) edge (v5)
		(v5) edge (v3)
		(v3) edge (v2)
		;
	\end{scope}

	\draw [to reversed-to](7,-0.3) -- (8,-0.3);

	\begin{scope}[xshift=9cm]
		\draw[gray!50](0,0.7) circle (0.2);
		\draw[gray!50](0,0) circle (0.2);
		\draw [rounded corners,gray!50] (-0.8,-0.55)--(0.8,-0.55)--(0,-1.7)--cycle;
		
		\node at (0,-2) {$\beta_1$};
		
		\vertex (v1) at (0,0.7) {};\node at (0.3,0.7) {$1$};
		\vertex (v2) at (0,0) {};\node at (0.3,0) {$2$};
		\vertex (v3) at (-0.5,-0.7) {};\node at (-0.8,-0.7) {$3$};
		\vertex (v4) at (0.5,-0.7) {};\node at (0.8,-0.7) {$4$};
		\vertex (v5) at (0,-1.4) {};\node at (0.3,-1.4) {$5$};
		
		\path
		(v4) edge (v5)
		(v4) edge (v3)
		(v5) edge (v3)
		(v3) edge (v2)
		;
	\end{scope}
	\end{tikzpicture}
	\caption{The $2$-simplex $\alpha_0$ is paired off with the $3$-simplex $\beta_0$ in $\V_3\setminus \V_{12}$, and the $2$-simplex $\alpha_1$ is paired off with the $3$-simplex $\beta_1$ in $\V_4\setminus \V_3$.}\label{counter}
\end{figure}
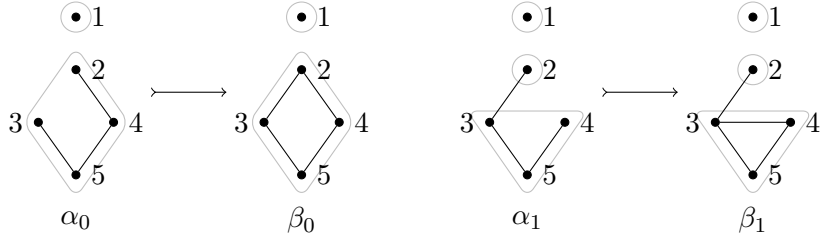
	From the construction of $\V$ on $\N_5$, it follows that $\alpha_0$ and $\beta_0$ are first paired off in $\V_3$, but $\alpha_1$ is not the head of an arrow, and $\alpha_1$ and $\beta_1$ are first paired off in $\V_4$.
\end{counterx}

\begin{counterx}
	It may also happen that $\alpha_0$ and $\beta_0$ are first paired in $\V_i$ for some $i \ge 3$, but $\alpha_1$ doesn't get paired \emph{at all} (and thus $\alpha_1$ is a critical simplex in $\V$). In that case the $\V$-path doesn't extend beyond $\alpha_1$, but not because $\alpha_1$ is the head of an arrow in $\V_3$.	A specific example for $n=4$ is in Figure~\ref{counter2}.
	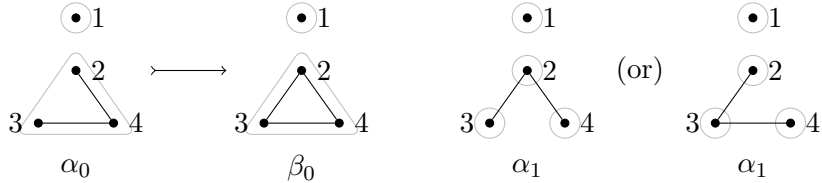
\begin{figure}[!ht]
		\centering
		\begin{tikzpicture}
			\begin{scope}			
				\draw[gray!50](0,0.7) circle (0.2);
				\draw [rounded corners,gray!50] (0,0.3)--(-0.8,-0.85)--(0.8,-0.85)--cycle;
				
				\node at (0,-1.3) {$\alpha_0$};
				
				\vertex (v1) at (0,0.7) {};\node at (0.3,0.7) {$1$};
				\vertex (v2) at (0,0) {};\node at (0.3,0) {$2$};
				\vertex (v3) at (-0.5,-0.7) {};\node at (-0.8,-0.7) {$3$};
				\vertex (v4) at (0.5,-0.7) {};\node at (0.8,-0.7) {$4$};
				
				\path
				(v2) edge (v4)
				(v4) edge (v3)
				;
			\end{scope}
			
			\draw [to reversed-to](1,0) -- (2,-0);
			
			\begin{scope}[xshift=3cm]
				\draw[gray!50](0,0.7) circle (0.2);
				\draw [rounded corners,gray!50] (0,0.3)--(-0.8,-0.85)--(0.8,-0.85)--cycle;
				
				\node at (0,-1.3) {$\beta_0$};
				
				\vertex (v1) at (0,0.7) {};\node at (0.3,0.7) {$1$};
				\vertex (v2) at (0,0) {};\node at (0.3,0) {$2$};
				\vertex (v3) at (-0.5,-0.7) {};\node at (-0.8,-0.7) {$3$};
				\vertex (v4) at (0.5,-0.7) {};\node at (0.8,-0.7) {$4$};
				
				\path
				(v2) edge (v4)
				(v4) edge (v3)
				(v3) edge (v2)
				;
			\end{scope}
			
			\begin{scope}[xshift=6cm]
				\draw[gray!50](0,0.7) circle (0.2);
				\draw[gray!50](0,0) circle (0.2);
				\draw[gray!50](-0.5,-0.7) circle (0.2);
				\draw[gray!50](0.5,-0.7) circle (0.2);
				
				\node at (0,-1.3) {$\alpha_1$};
				
				\vertex (v1) at (0,0.7) {};\node at (0.3,0.7) {$1$};
				\vertex (v2) at (0,0) {};\node at (0.3,0) {$2$};
				\vertex (v3) at (-0.5,-0.7) {};\node at (-0.8,-0.7) {$3$};
				\vertex (v4) at (0.5,-0.7) {};\node at (0.8,-0.7) {$4$};
				
				\path
				(v2) edge (v4)
				(v3) edge (v2)
				;
			\end{scope}
			
			\node at (7.5,0) {(or)};
			
			\begin{scope}[xshift=9cm]
				\draw[gray!50](0,0.7) circle (0.2);
				\draw[gray!50](0,0) circle (0.2);
				\draw[gray!50](-0.5,-0.7) circle (0.2);
				\draw[gray!50](0.5,-0.7) circle (0.2);
				
				\node at (0,-1.3) {$\alpha_1$};
				
				\vertex (v1) at (0,0.7) {};\node at (0.3,0.7) {$1$};
				\vertex (v2) at (0,0) {};\node at (0.3,0) {$2$};
				\vertex (v3) at (-0.5,-0.7) {};\node at (-0.8,-0.7) {$3$};
				\vertex (v4) at (0.5,-0.7) {};\node at (0.8,-0.7) {$4$};
				
				\path
				(v4) edge (v3)
				(v3) edge (v2)
				;
			\end{scope}
		\end{tikzpicture}
		\caption{The $1$-simplex $\alpha_0$ is paired off with the $2$-simplex $\beta_0$ in $\V_3\setminus \V_{12}$. Both choices of $\alpha_1$ are unpaired in $\V$.}\label{counter2}
	\end{figure}
	Here $\alpha_0$ and $\beta_0$ are first paired off in $\V_3$. There are two choices for $\alpha_1$, but both of them are unpaired (i.e., critical) in $\V$.
\end{counterx}

	\section{Alternative proof of $\V$ being a gradient vector field}\label{alt-proof}
	We now provide an alternative justification for the claim that $\V$ is indeed a gradient vector field on $\N_n$. As mentioned before, $\V_{12}$ is a gradient vector field on $\N_n$. Let, if possible, $\alpha_0, \beta_0, \alpha_1, \beta_1, \ldots, \alpha_r, \beta_r, \alpha_{r+1}=\alpha_0$ be a non-trivial closed $\V$-path, denoted by $\gamma$. Also, suppose $\alpha_0$ and $\beta_0$ are first paired off in $\V_k$, with $k \ge 3$. It follows that $\alpha_0$ has two components, and $\beta_0=\alpha_0+(j,k)$, for some $j<k$.
	
	Let, if possible, the set $\{1,\ldots,r\}$ contain a least element, say $\ell$, such that $\alpha_\ell$ has more than two components. Then it follows that the vertices $1$ and $2$ are in two different components of $\alpha_\ell$, and $\beta_\ell=\alpha_\ell+(1,2)$. This implies, as $\alpha_{\ell+1}$ is different from $\alpha_\ell$, the graph $\alpha_{\ell+1}$ contains the edge $(1,2)$, and thus $\alpha_{\ell+1}$ is the head of an arrow (in $\V_{12}$). Therefore, $\gamma$ cannot be continued further, a contradiction. So suppose each $\alpha_i$ (and thus each $\beta_i$) has exactly two components.
	
\begin{figure}[!ht]
		\centering
		\begin{tikzpicture}
			\node at (1,-2) {$\alpha_0$};
			\begin{scope}
				\draw [rounded corners,dashed,blue!90] (0,1)--(-0.9,0)--(0,-1)--(0.9,0)--cycle;

				\vertex (v1) at (0,0.7) {};\node at (0.3,0.7) {$1$};
				\vertex (v3) at (-0.7,0) {};\node at (-1,0) {$3$};
				\vertex (v5) at (0.7,0) {};\node at (0.4,0) {$5$};
				\vertex (v6) at (0,-0.7) {};\node at (0.3,-0.7) {$6$};
				\vertex (v9) at (-0.7,-1) {};\node at (-0.7,-1.3) {$9$};
				\vertex (v11) at (-1.4,-0.7) {};\node at (-1.7,-0.7) {$11$};
				
				\path
				(v1) edge (v3)
				(v1) edge (v5)
				(v3) edge (v6)
				(v3) edge (v9)
				(v3) edge (v11)
				(v9) edge (v11)
				;
			\end{scope}
		
		\begin{scope}[xshift=2cm]
			\draw [dashed,blue!90] (0,0.4) ellipse (0.3cm and 0.6cm);
			\draw [rounded corners,red!70] (0,0.3)--(-0.7,-0.7)--(0,-1.3)--(0.7,-0.7)--cycle;
			
			\vertex (v2) at (0,0.7) {};\node at (0.3,0.7) {$2$};
			\vertex (v4) at (0,0) {};\node at (0.3,0) {$4$};
			\vertex (v7) at (-0.5,-0.7) {};\node at (-0.8,-0.7) {$7$};
			\vertex (v8) at (0.5,-0.7) {};\node at (0.8,-0.7) {$8$};
			\vertex (v10) at (0,-1) {};\node at (0,-1.3) {$10$};
			
			\path
			(v2) edge (v4)
			(v4) edge (v8)
			(v4) edge (v10)
			(v8) edge (v10)
			(v10) edge (v7)
			;
		\end{scope}
		
		\end{tikzpicture}
		\caption{The graph $\alpha_0$ ($\in \N_{11}$) first gets paired off with $\alpha_0+(4,7)$ in $\V_7$, and thus, $k=7$, $j=4$. Here, $F$ is the subgraph induced by $\{1,\ldots,6\}$ (bounded by dashed (blue) lines), and $H_0$ is the subgraph induced by $\{4,7,8,10\}$ (bounded by solid (red) lines).}\label{FH}
	\end{figure}
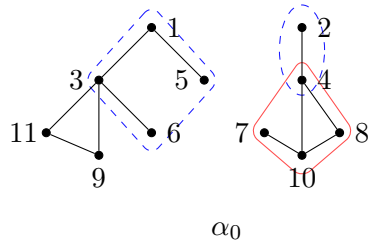

	Let $F$ be the subgraph of $\alpha_0$, induced by the vertex set $\{1,\ldots,k-1\}$ (see Figure~\ref{FH}). From the construction of $\V$, it follows that $F$ is a forest (on the vertex set $\{1,\ldots,k-1\}$) with two components such that the vertices $1$ and $2$ are in different components. Moreover, both components of $F$ have the property that the vertex labels are increasing along every ray starting from the vertex $1$ or the vertex $2$. Also, since each $\alpha_i$ has two components and $F$ is a disjoint union of two trees, it follows that the subgraph of each $\alpha_i$, induced by the vertex set $\{1,\ldots,k-1\}$, is the same as $F$. This also implies that no $\alpha_i$ is paired off in $\V_{k-1}$ (we abuse the notation when $k=3$; by $\V_2$ we mean $\V_{12}$).
	
	Let $H_i$ be the component of the graph $\alpha_i - E(F)$ (i.e., the graph obtained from $\alpha_i$ after deleting the edges in $F$) that contains the vertex $k$ (see Figure~\ref{FH}). We note that out of the vertices $1,\ldots,k-1$, only $j$ is a vertex of $H_0$, and consequently, $j$ is a vertex of each $H_i$. Since $\gamma$ is a non-trivial closed $\V$-path, and $\alpha_1$ contains the edge $(j,k)$ (but $\alpha_0$ doesn't), the set $\{2,\ldots,r+1\}$ contains a least element, say $\ell$, such that $\alpha_\ell$ doesn't contain the edge $(j,k)$ (i.e., $\alpha_\ell=\beta_{\ell-1}-(j,k)$). As the connected graph $H_\ell$ contains the vertex $j$, and $\alpha_\ell$ is unpaired in $\V_{k-1}$, the simplex $\alpha_\ell$ must be paired off with $\alpha_\ell+(j,k)=\beta_{\ell-1}$, which leads to a contradiction. This concludes the proof.
	
	\section{Concluding remarks}\label{concl}
	Here we remark that although discrete Morse theoretic techniques are purely combinatorial and
	generally computationally efficient in nature, one needs to be meticulous while applying them. While exploring the structure of a complex via these techniques and performing various computations, we are usually required to start with a sufficiently ``good" gradient vector field. One primary challenge often lies in showing that a constructed discrete vector field is a gradient one (i.e., justifying its ``acyclicity"). Our proof makes use of the following observation which is useful for verifying the acyclicity of a constructed discrete vector field in general.
	\begin{obs}\label{obs}
		Let $\V$ be a discrete vector field on a complex and $\gamma :$ $\alpha_0,$ $\beta_0,$ $\alpha_1,$ $\ldots,$ $\alpha_k,$ $\beta_k,$ $\alpha_{k+1}$ be a $\V$-path such that $\beta_0 = \alpha_0 \cup \{x\}$ (with $x \notin \alpha_0$) and $\alpha_1 = \beta_0 \setminus \{y\}$ (with $y \in \alpha_0$). If $\gamma$ is a closed $\V$-path, then
		\begin{enumerate}[(i)]
			\item the set $\{2,\ldots,k\}$ contains an (least) element $i$ such that $\alpha_i$ doesn't contain $x$, and
			\item the set $\{3,\ldots,k+1\}$ contains an (least) element $j$ such that $y \in \alpha_j$.
		\end{enumerate}
		In other words, if the assumption that there is an (least) $i \in \{2,\ldots,k\}$ such that $x \notin \alpha_i$ (or, there is a (least) $j \in \{3,\ldots,k+1\}$ such that $y \in \alpha_j$), leads to a contradiction for an arbitrary $\gamma$, then there is no closed $\V$-path. Consequently, $\V$ is a gradient vector field.
	\end{obs}
	
	\section*{Acknowledgements}
	The authors would like to express their sincere gratitude to Prof.\ Goutam Mukherjee and Sajal Mukherjee for their encouragement and support throughout this work. The first author would also like to thank Apratim Chakraborty and Kuldeep Saha, in addition to S.~Mukherjee, for their insightful discussions on discrete Morse theory during previous collaborations. These discussions significantly deepened our understanding of this subject area.

\end{document}